\documentclass{amsart}
\usepackage{amsmath,amssymb,amsfonts}
\usepackage[latin1]{inputenc}
\usepackage[english]{babel}
\usepackage{graphicx,xcolor} 
\graphicspath{ {./Figures/} }
\usepackage{float} 
\usepackage[figuresright]{rotating}
\usepackage{color}

\usepackage{hyperref}

\newcommand{\RR}{{\mathbb R}}

\newcommand{\EEE}{{\rm I\kern-2pt E}}
\newcommand{\NN}{{\mathbb N}}
\newcommand{\PNS}{P^{j}_s}

\newcommand{\PP}{{\rm I\kern-2pt P}}

\newcommand{\PNT}{P^{j}_t}

\newcommand{\ANT}{\alpha^{j}_t}

\newcommand{\Hmin}{{H^{min}_f}}
\newcommand{\ome}{\omega}
\newcommand{\Ome}{\Omega}

\newcommand{\la}{{\lambda}}
\newcommand{\La}{{\Lambda}}

\newcommand{\ep}{{\varepsilon}}
\newcommand{\al}{{\alpha}}
\newcommand{\be}{\beta}

\newcommand{\ZZ}{{\mathbb Z}}
 
\newcommand{\cli}{c_{j,k}}

\newcommand{\BE}{\begin{equation}}
\newcommand{\EE}{\end{equation}}

\newtheorem{theo}{Theorem} 
\newtheorem{prop}{Proposition}[section]
\newtheorem{coro}[prop]{Corollary} 
\newtheorem{lemm}[prop]{Lemma}
\newtheorem{defi}{Definition}

\newcommand{\BP}{\begin{prop}}
\newcommand{\EP}{\end{prop}}
\newcommand{\BL}{\begin{lemm}}
\newcommand{\EL}{\end{lemm}}
\newcommand{\BD}{\begin{defi}}
\newcommand{\ED}{\end{defi}}
\newcommand{\BT}{\begin{theo}}
\newcommand{\ET}{\end{theo}}
\newcommand{\BC}{\begin{coro}}
\newcommand{\EC}{\end{coro}}

\begin{document}

\title{\boldmath Multifractal analysis based on $p$-exponents and lacunarity exponents}

%


\author[P. Abry]{Patrice Abry}
\address{%
Patrice Abry, Signal, Systems and Physics, Physics Dept., CNRS UMR 5672, Ecole Normale Sup\'erieure de Lyon, Lyon, France} 
\email{patrice.abry@ens-lyon.fr}

\author[S. Jaffard]{St\'ephane Jaffard}
\address{%
St\'ephane Jaffard, Universit\'e Paris Est, Laboratoire d'Analyse et de Math\'ematiques Appliqu\'ees, CNRS UMR 8050, UPEC,  Cr\'eteil, France}
\email{jaffard@u-pec.fr}

\author[R. Leonarduzzi]{Roberto Leonarduzzi}
\address{%
Roberto Leonarduzzi, Signal, Systems and Physics, Physics Dept., CNRS UMR 5672, Ecole Normale Sup\'erieure de Lyon, Lyon, France}
\email{roberto.leonarduzzi@ens-lyon.fr}

\author[C. Melot]{Clothilde Melot}
\address{%
Clothilde Melot, I2M,  CNRS UMR 7373, Universit\'e Aix-Marseille, Marseille,  France}  
\email{melot@cmi.univ-mrs.fr}

\author[H. Wendt]{Herwig Wendt}
\address{%
Herwig Wendt, IRIT, CNRS UMR 5505, University of Toulouse,  France}
\email{herwig.wendt@irit.fr}


\begin{abstract} 

Many examples of signals and images cannot be modeled by locally bounded functions, so that the standard multifractal analysis, based on the H\"older exponent,  is not feasible. We  present a multifractal analysis based on another quantity, the $p$-exponent, which can take arbitrarily large negative values. We investigate some mathematical properties of this exponent, and show how it allows us  to model the idea of ``lacunarity'' of a singularity at  a point. We finally adapt the wavelet based  multifractal analysis  in this setting, and we give applications to a simple mathematical model of multifractal processes: Lacunary wavelet series.
\end{abstract}

\maketitle


\section{Introduction}

The origin of fractal geometry   can be traced back to     the quest for non-smooth functions, rising from  a key question that motivated a large part of the progresses in analysis during the nineteenth century: Does a continuous function necessarily have points of differentiability? A negative answer to this question was supplied by Weierstrass when he built his famous counterexamples, now referred  to as the { \em Weierstrass functions}

 \BE \label{weier1}  \mathcal{W}_{a, b} (x) = \sum_{n =0}^{+ \infty}  a^n {\cos ( b^n \pi x)} \hspace{15mm}  \EE
 where $0 < a <1$, $b$ was an odd integer  and $ab  > 1 + 3 \pi /2$. The fact that they are continuous and nowhere differentiable was  later sharpened by Hardy   in a way which requires
  the notion of { \em pointwise  H\"older regularity}, which is the most commonly used notion of pointwise regularity in the function setting. We assume in the following that the functions or distributions we consider are defined on $\RR$. However,  most results that we will  investigate extend to several variables.

 \BD  \label{defholdreg}  Let $f: \; \RR  \rightarrow  \RR $  be a { locally   bounded function}, $x_0 \in   \RR$ and  let $\gamma \geq 0$; $f $ belongs to $ C^\gamma (x_0)$ if there exist $C>0$, $R >0$  and a polynomial  $P$ of degree less than  $\gamma$ such that:  
\BE \label{equ:tpe}  \mbox{ for a.e. }\;\;   x \;\; \mbox{ such that  } \;\; |x-x_0| \leq R,  \qquad  |f(x) -P(x-x_0) | \leq C | x-x_0|^\gamma . \EE  
The H\"older exponent of    $f$ at $x_0
$ is  
\BE \label{equ:he}
  h_f  (x_0) =\sup \left\{ \gamma : \;\; f \;\; \mbox{ is  }   \;\;   C^{\gamma
} (x_0) \right\} . \EE
\ED

The H\"older exponent of $\mathcal{W}_{a, b}$ is a constant function, which is equal to $H= -\log a / \log b $ at every point (see e.g.  \cite{Jaffard2004} for a simple, wavelet-based proof); since $H <1$ we thus recover the fact that $\mathcal{W}_{a, b}$ is nowhere differentiable, but the sharper notion of 
H\"older exponent allows us to draw a difference between each of  the   Weierstrass functions, and classify them  using a regularity parameter that takes values in $\RR^+$.  
 The  graphs of Weierstrass functions  supply important examples of fractal sets that  still motivate research  (the determination of their Hausdorff dimensions  remains   partly open, see \cite{Baransky}). 
 In applications, such fractal characteristics have been used for classification purposes. For instance, an unorthodox use was the discrimination between Jackson Pollock's original paintings and fakes using the box dimension of the graph supplied by the pixel by pixel values of a high resolution photograph of the painting,  see  \cite{rcew2008}. 
 
The status of everywhere irregular functions was, for a long time, only the one of academic counter-examples, such as the Weierstrass functions. This situation changed when stochastic processes like  Brownian motion  (whose H\"older exponent is $H=1/2$ everywhere) started to play  a key role in the modeling of physical phenomena. Nowadays, experimentally acquired signals that are everywhere irregular are prevalent in a multitude of applications, so that the classification and modeling of such data  has become a key problem. 
 However,   the use of a single parameter (e.g. the box dimension of the graph) is too reductive as a classification tool in many situations that are met in applications.  
 This explains the success of multifractal analysis, which is a way to associate a whole collection of fractal-based  parameters to a function.
 Its purpose  is twofold:  on the mathematical side, it allows one to determine the size of the sets of points where a function has a given H\"older exponent; on the signal processing side,   it yields  new collections of parameters  associated to the considered signal and which can be used for  classification, model selection, or for parameter selection inside a parametric setting.  The main advances in the subject came from  a better  understanding of the interactions between these two motivations, e.g., see \cite{MandMemor} and references therein for recent review papers. 

Despite the fact that multifractal analysis has traditionally been based on the Hölder exponent, it is not the only characterization of pointwise regularity that can be used. Therefore, our  goal in the present contribution is to analyze alternative pointwise exponents and the information they provide. 

 In Section \ref{secpointexp}  we  review  the possible pointwise exponents  of functions, and explain in which context  each  can be used.
 
 In Section \ref{fracexp} we  focus on the $p$-exponent, derive some of its properties, and  investigate what information it yields concerning the lacunarity of the local behavior of the function near a singularity.
 
 In Section \ref{seclacu} we recall the derivation of the multifractal formalism  and give applications to a simple model of a random process which displays multifractal behavior: Lacunary wavelet series.  
 
  We conclude with remarks  on the relationship between the existence of $p$-exponents and the sparsity of the wavelet expansion. \\
  
  This paper partly reviews elements on the  $p$-exponent which are scattered in the literature, see e.g. \cite{Abel,CalZyg,JaffCies,JaffToul,PART1}. New material starts with the introduction and analysis of the lacunarity exponent in Section \ref{sec:lacunary}, the analysis of thin chirps in Section \ref{thinch}, and all following sections, except for the brief reminder on the multifractal formalism in Section \ref{derivfp}. 
 
 
 \section{Pointwise exponents  } 
 
 
 \label{secpointexp} 
 


In this section, unless otherwise specified, we assume that $f \in L^1_{loc} (\RR )$. 
An important remark concerning the definition of pointwise H\"older regularity  is 
 that if (\ref{equ:tpe}) holds (even for $\gamma <0$), then $f$ is  bounded in any annulus  $0 < r \leq | x-x_0| \leq R $. It follows that, if an estimate such as (\ref{equ:tpe}) holds for all $x_0$, then $f$ will be locally bounded, except perhaps at isolated points. 
For this reason, one usually assumes that  the considered  function $f$ 
 is  (everywhere)  locally bounded. It follows that (\ref{equ:tpe}) holds for $\gamma =0$ so that  the H\"older exponent is always nonnegative. 
 
 \subsection{Uniform H\"older regularity} 
 
 An important issue therefore is to determine if the  regularity assumption $f \in L^\infty_{loc}$ is satisfied for real life data.  This can be done in practice by first determining their  { \bf uniform H\"older exponent},   which is defined as follows. 
 

   Recall that  Lipschitz spaces $C^s (\RR ) $ are defined for $0 < s <1$  by 
\[ f \in L^\infty \quad \mbox{ and} \quad  \exists C ,  \;\;  \;\; \forall x,y,  \qquad | f(x)-f(y) | \leq C |x-y|^s . \]
If $s >1$, they are then defined by recursion on $[s]$ by the condition:
$ f \in C^s (\RR )$ if $f\in L^\infty$ and  if  its  derivative $f'$ (taken in the sense of distributions)  belongs to 
$C^{s-1} (\RR)$. If $s <0$, then the $C^s$ spaces are composed of distributions,  also defined by recursion on $[s]$  as follows:
$ f \in C^s  (\RR )$ if $f$ is a derivative (in the sense of distributions) of a function 
$g \in C^{s+1} (\RR )$. 
We thus obtain a definition of   the $C^s$ spaces  for any $s \notin \ZZ$ (see \cite{Mey90I} for   $s\in \ZZ$,  which we will however not need  to consider   in the following).  
A distribution $f$ belongs to $C^{s }_{loc}$ if $f \varphi \in C^s$ for every $C^\infty$ compactly supported function $\varphi$.

\BD The  uniform H\"older exponent  of  a tempered distribution $f$ is 
\BE \label{caracbeswav2hol}  \Hmin = \sup \{ s : \; f \in C^{s }_{loc} (\RR) \}   . \EE
\ED
This definition does not make any a priori assumption on $f$: 
The uniform H\"older exponent  is defined for any tempered distribution, and it can be positive or negative. 
More precisely:   
 \begin{itemize}
\item  If {  $\Hmin > 0$}, then $f $ is  a locally bounded function,
\item  if {$\Hmin <0$}, then $f$ is not a locally bounded function.
\end{itemize}
 
In practice,  this exponent is determined through the help of the wavelet coefficients of  $f$. 
 By definition, an { \bf orthonormal wavelet basis} is generated  by a couple of functions
$(\varphi, \; 
\psi^{}_{})$, which, in our case, will   either be in the Schwartz class, or   smooth and  compactly supported (in that case, wavelets are assumed to be smoother than the regularity exponent of the considered space). The  functions $\varphi (x
-k), \;\;\ k\in \ZZ, $ together with 
$  2^{j/2} \psi (2^jx
-k), \;\;\; j\geq 0, \; k\in \ZZ ,$ form an orthonormal basis of $L^2( \RR
 ) $. Thus any function $f\in L^2(\RR )$ can be written 
\[ f(x) =\sum_{k} c_k  \;  \varphi(x-k)+  \sum_{j\geq 0 } \sum_{k \in \ZZ} \cli \;  \psi(2^jx-k), \]
where the wavelet coefficients of $f$ are given by 
\BE \label{defcoef}  c_k =\int  \varphi (t-k) f(t) dt  \quad \mbox{ and} \quad  \cli =2^{j}  \int  \psi (2^jt-k) f(t) dt .\EE
An important remark is that these formulas also hold in many different functional settings (such as the Besov or Sobolev spaces of positive or negative regularity),  provided that the picked wavelets are smooth enough (and that the integrals (\ref{defcoef}) are understood as duality products).

Instead of using  the  indices $(j,k)$, we will often use dyadic intervals: Let  
\BE \label{deflamb}   \lambda\; (= \lambda (j,k)) \; = 
  \left[ \frac{k}{2^j} , \frac{k+1}{2^j} \right) \EE
  and, accordingly: $c_{\lambda} =\cli$ 
 and   $\psi_{\lambda} (x)= \psi (2^jx-k) $. 
  Indexing by dyadic intervals will be   useful  in the sequel because  the interval $\lambda$ indicates the localization of the corresponding wavelet:  
 When the wavelets are compactly supported, then, 
$ \exists C >0 $ such that  when 
$supp (\psi ) \subset [ -C/2, C/2 ]$, then
$supp (\psi_\lambda  )  \subset  2C \lambda.$  \\


 In practice, $ \Hmin $ can be derived  directly from the wavelet coefficients of $f$ through a simple regression in a log-log plot; 
indeed, it follows from the wavelet characterization of the spaces $C^s$, see \cite{Mey90I}, that:
\BE 
\label{caracbeswav3hol}  
\Hmin =\liminf_{j \rightarrow + \infty} \;\;  
 \;\;  \frac{ \log 
\left(  \displaystyle  \sup_{k}  | \cli |   \right) }{\log (2^{-j})}. 
\EE
This estimation procedure has been studied in more detail in \cite{Bergou}. Three examples of its numerical application to real-world functions are provided in Figure \ref{fig:rwimage}.

A multifractal analysis based on the H\"older exponent can only be performed if $f$ is locally bounded. 
A way to determine if this is the case consists in first checking if $\Hmin >0$. 
 This quantity  is perfectly well-defined for mathematical functions  or stochastic processes; e.g. for Brownian motion,  $\Hmin =1/2$, and for Gaussian white noise,  $\Hmin =-1/2$. However the situation may seem less clear for experimental signals; indeed any data acquisition device yields a finite set of locally averaged quantities, and one  may argue that  such a finite collection of data (which, by construction, is bounded) can indeed be modeled by a locally bounded function. This argument can only be turned by revisiting the way that \eqref{caracbeswav3hol} is computed in practice: Estimation is performed through a linear regression in log-log coordinates { \bf on the range of scales available in the data} and $\Hmin$ can indeed be found  negative for a finite collection of data. At the modeling level, this means that a mathematical model which would display the same linear behavior in log-log coordinates { \bf at all scales} would satisfy $\Hmin <0$.

 The quantity $\Hmin$ can be found either positive or negative depending on the nature of the application.
For instance, velocity turbulence data  and price time series in finance  are found to always have $ \Hmin>0$, while aggregated count Internet traffic time series always have $ \Hmin<0$. 
For biomedical applications (cf. e.g., fetal heart rate variability) as well as for image processing, $ \Hmin $ can commonly be  found either positive or negative (see Figure \ref{fig:rwimage}) \cite{AWJH,MandMemor,Bergou,Zuhai,WENDT:2009:C}. 
This raises the problem of using other pointwise regularity exponents that would not require the assumption that the data are locally bounded. We now introduce such exponents.




 \subsection{\boldmath The $p$-exponent for $p \geq 1$ }
 
 The introduction of $p$-exponents is motivated by the necessity of introducing  regularity exponents  that could be defined even when $\Hmin$ is found to be negative;    $T^p_\alpha (x_0)$  regularity,  introduced  by A. Calder\'on and A. Zygmund in \cite{CalZyg},   has the advantage of only making the assumption that $f$   locally belongs to  $ L^p (\RR)$. 

 \BD Let $p \geq 1$ and assume that $f\in  L^p_{loc} (\RR)$.  
Let  $\alpha \in \RR$;  
the  function $ f $   belongs to $T^p_\alpha (x_0)$  if there exists $ C$  and a polynomial $P_{x_0}$  of degree  less than $\alpha $  such that, for $r$ small enough,  
\begin{equation}  
\label{pexpa} 
 \left( \frac{1}{2r} \int_{x_0 -r}^{x_0 +r} | f(x) -P_{x_0}(x)|^p dx  \right)^{1/p} \leq C r^\alpha.
\end{equation}
\ED

Note that the { \bf Taylor polynomial} $P_{x_0}$  of $f$ at $x_0$ might depend on $p$. However, one can check that only its degree does (because the best possible $\al$ that one can pick in (\ref{pexpa}) depends on $p$ so that its integer part  may vary with $p$, see \cite{Abel}). Therefore  we introduce no such dependency in the notation, which will lead to no ambiguity afterwards. 

The { \bf $p$-exponent}   of  $f$ at $x_0$ is defined as 
\begin{equation}  
\label{equ-pexp} 
h_f^p(x_0) = \sup \{ \alpha : f \in T^p_\alpha (x_0)\}. 
\end{equation}

The condition  that $f$  locally belongs to  $ L^p (\RR)$ implies that (\ref{pexpa}) holds for $\alpha = -1/p$, so that $h_f^p (x_0)  \geq  -1/p$. 
\\

We will consider in the following ``archetypical'' pointwise singularities, which are simple toy-examples of singularities with a specific behavior at  a point. 
They will  illustrate the new notions we consider and  they will also supply  benchmarks on which we can compute exactly what these new notions allow us  to quantify. These toy-examples  will be a    test for the adequacy between these  mathematical notions and the intuitive  behavior that we expect to quantify. The first (and most simple)  ``archetypical'' pointwise singularities are the { \bf  cusp singularities}.  

Let $\al \in \RR -2 \NN$ be such that $\al > -1$. The { \bf cusp of order $\al$ at $0$}  is the function 
\BE \label{cusp} { \mathcal C}_{\al} (x) = | x|^\al .  \EE  
The case $\al \in 2 \NN$ is excluded because it leads to a $C^\infty$ function. However, if $\al = 2n $, one can pick 
\[  { \mathcal C}_{2n} (x) =  x | x|^{2n-1} , \]  
in order to cover this case  also. 

If $\al \geq 0$, then the cusp ${ \mathcal C}_{\al}$ is locally bounded and  its H\"older exponent  at $0$   is well-defined and takes the value $\al$.  If   $\al  > -1/p$, then   its  $p$-exponent  at $0$ is well-defined and also takes the value $\al$, as in the H\"older case. (Condition $\al  > -1/p$ is  necessary and sufficient to ensure that  ${ \mathcal C}_{\al}$ locally belongs to $L^p$.) Examples for cusps with several different values of $\alpha$ are plotted in Figure \ref{fig:cusp}.

If $f \in L^p_{loc}$ in a neighborhood of $x_0$ for a $p \geq 1$, let us define the { \bf critical Lebesgue index} of $f$ at $x_0 $ by 
 \BE  \label{critleb} p_0 (f) = \sup \{ p: f \in  L^{p}_{loc}  (\RR)  \mbox{ in a neighborhood of   $x_0$}  \} .  \EE  
 The importance of this exponent comes from the fact that it tells in practice for which values of $p$ a $p$-exponent based multifractal analysis can be performed. Therefore, its numerical determination is an important prerequisite that should not be bypassed in applications. In Section \ref{subsecplone} we will    extend  the definition of $p_0 (f)$  to situations where $f \notin L^1_{loc}$ and show how it can be derived from another quantity, the wavelet scaling function, which can be effectively computed on real-life data.
 
 \subsection{The lacunarity exponent} 
 \label{sec:lacunary}
 
 The $p$-exponent at $x_0$ is  defined on the  interval  $[1, p_0 (f)  ]$ or $[1, p_0(f)  )$;  
when the $p$-exponent does not depend on $p$ on this interval, we will say   that $f$ has a { \bf\boldmath  $p$-invariant singularity} at $x_0$.  Thus, cusps are 
$p$-invariant singularities.

This first example raises the following question: Is  the notion of $p$-exponent only relevant as an  extension of  the  H\"older exponent to non-locally bounded functions?   Or  can it take different values with $p$, even for bounded functions?  And, if such is the case, how can one characterize the additional information thus supplied?  In order to answer this question, we introduce a second type of archetypical singularities, the { \bf lacunary singularities},
which will show that the $p$-exponent may be non-constant.  We first need to recall  the  geometrical notion of { \bf accessibility exponent}  which quantifies the { \bf lacunarity} of a set at a point, see \cite{JaffMel}. We denote by $\mathcal{M}(A)$ the Lebesgue measure of a set $A$. 

\begin{defi}  \label{def:access} Let $\Ome\subset \RR $.  
 A point $x_0$ of the boundary of $\Ome$ is    
$\al$-accessible if  there exist $C>0$  and $r_0>0$ such that $\forall r \leq
r_0$, 
\begin{equation} 
\label{eq1}  \mathcal{M} \left( \Ome\cap B(x_0,r) \right) \leq C r^{\al+1}.
\end{equation} 

The supremum of all values of $\al$ such that (\ref{eq1}) holds is called
the  {accessibility exponent of $\Ome$  at $x_0$}.  We will denote it by ${ \mathcal E}_{x_0} ( \Omega )$. 
 \end{defi}

Note that ${ \mathcal E}_{x_0} ( \Omega ) $ is always nonnegative. 
If it is strictly positive, then $\Omega$ is lacunary at $x_0$.
The accessibility exponent  
  supplies a way to estimate, through a log-log plot regression,  the ``size'' of the part of $\Omega$  which is contained in arbitrarily small neighborhoods of $x_0$.  The following sets illustrate this notion.

  Let  $\ome$ and $ \gamma$ be such that $0 < \gamma \leq \ome$;   the  set $U_{\ome, \gamma} $  is defined as follows. Let 
 \BE \label{singlacex}  I^j_{\ome, \gamma} =    [2^{-\ome j} , 2^{-\ome j} + 2^{-\gamma j}]; \quad \mbox{ then}  \quad U_{\ome, \gamma}  = \bigcup_{j \geq 0}   I^j_{\ome, \gamma} . \EE
 Clearly, at the origin,  
 \BE \label{lac1} { \mathcal E}_{0} ( U_{\ome, \gamma} ) = \frac{\gamma}{\ome} -1. \EE

 We now construct  univariate functions $F_{ \al, \gamma}: \RR \rightarrow \RR$ which  permit us to better understand the conditions under which $p$-exponents will differ. These functions will have a lacunary support in the sense of Definition \ref{def:access}. 

Let $\psi$ be the Haar wavelet: $\psi =1_{[0, 1/2)} - 1_{[1/2, 1)}$  and \[ \theta(x) = \psi (2x) - \psi (2x-1) \] (so that $\theta$ has the same support as $\psi$ but its two first moments vanish). 

\BD 
Let $\alpha \in \RR  $ and $\gamma  >1$. The { lacunary comb  }   $F^\al_{ \ome, \gamma}$ is the function  
\BE \label{defaog}  F^\al_{ \ome, \gamma} (x) = \sum_{j=1}^\infty 2^{-\al j} \theta \left( 2^{\gamma j} (x-2^{-\ome j})  \right) .  \EE
\ED

\begin{figure}[tb]
\includegraphics[width=\linewidth]{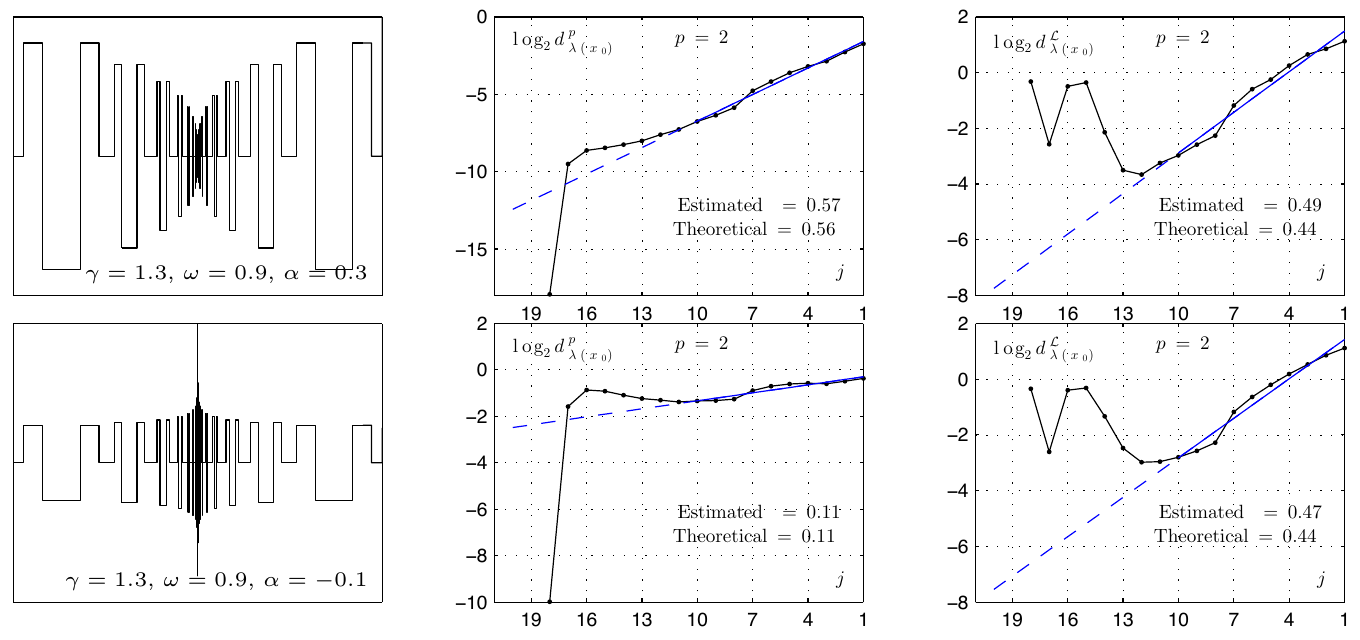}
\caption{\label{fig:Lcomb}{\bf Lacunary combs} with $p_0=+\infty$ (top row) and $p_0=9.1$ (bottom row): functions (left column) and estimation of $p$-exponents and lacunarity exponents (center and right column, respectively).}
\end{figure}

Note that its singularity is at $x_0 =0$. 
Numerical examples of lacunary combs are provided in Figure \ref{fig:Lcomb}.

Note  that the support of $F^\al_{ \ome, \gamma} $ is  $U_{\ome, \gamma}$  so that the  accessibility exponent  at 0  of this  support  is given by (\ref{lac1}). 
The function  $F^\al_{ \ome, \gamma} $  is locally bounded if and only if $\al \geq 0$. Assume that $\al <0$; then $F^\al_{ \ome, \gamma}  $ locally belongs to $L^p$ 
if and only if $\al > -\gamma /p$. 
 When such is the case, a straightforward computation yields that its $p$-exponent  at 0 is given by
 \BE \label{pexpcomb}  h_{F^\al_{ \ome, \gamma} }^p (x_0) = \frac{\al}{\ome} + \left( \frac{\gamma}{\ome } - 1 \right) \frac{1}{p}.\EE 
 In contradistinction with the cusp case, the $p$-exponent    of  $F^\al_{ \ome, \gamma} $  at 0 is not a constant function of $p$.  Let us see how the variations of the mapping $p \rightarrow h^p_f (x_0)$ are related  with the lacunarity of the support of $f$, in the particular case of  $F^\al_{ \ome, \gamma} $.  
 We note that this mapping is an affine function of the variable $q = 1/p$ (which, in this context, is  a more natural parameter than  $p$)  and that  the accessibility exponent of the support of $F^\al_{ \ome, \gamma} $ can be recovered by a derivative of this mapping with respect to  $q$. 
 The next question is to determine the value of $q$ at which this derivative should be taken.  This toy-example is too simple to   give a clue since any value of $q$ would lead to the same value for the derivative. We want to find if there is a more natural one, which would lead to a { \em canonical  } definition for the lacunarity exponent.  It is possible to settle this point through the following  simple perturbation argument:   Consider a new singularity $F$
 that would be the sum of two functions $F_1= F^{\al_1}_{ \ome_1, \gamma_1}$ and $F_2= F^{\al_2}_{ \ome_2, \gamma_2}$ with 
 \BE \label{condalga}  0 < \al_1 < \al_2 \quad \mbox{  and }  \quad 
 \gamma_1 > \gamma_2. \EE
 The $p$-exponent of $F$ (now expressed in the $q$ variable, where $q = 1/p$) is given by
 \BE \label{deriv}  q\mapsto h_f^{\frac{1}{q}}(x_0)=\min \left[  \frac{\al_1}{\ome} + \left( \frac{\gamma_1}{\ome } - 1 \right) q\; , \;  \;    \frac{\al_2}{\ome} + \left( \frac{\gamma_2}{\ome } - 1 \right) q \right] . \EE  
 The formula for  the lacunarity exponent  should yield  the lacunarity  of the most irregular component of $F$; since $F \in L^\infty_{loc}$,  the H\"older exponent   is the natural way to measure this irregularity.
 In this respect,  the most irregular component  is  $F_1$; the lacunarity exponent should thus take the value   $\left( \frac{\gamma_1}{\ome } - 1 \right) $. But, since (\ref{condalga}) allows the shift  in slope of the function (\ref{deriv}) from $ \left( \frac{\gamma_1}{\ome } - 1 \right)$ to $\left( \frac{\gamma_2}{\ome } - 1 \right)$  to take place  at a $q$  arbitrarily close to $0$,  the only way to  obtain this desired result in any case is to pick the derivative of the mapping $q \rightarrow h^{1/q}_f (x_0)$  precisely at $q=0$.\\

 A similar perturbation   argument can be developed  if $p_0 (f) < \infty$  with the conclusion  that the derivative should be  estimated  at the smallest possible value of $q$, i.e. for 
 \[ q= q_0 (f) := \frac{1}{p_0 (f)}; \] 
 hence the following   definition of  the {\bf lacunarity exponent}.
 
 \BD   \label{deflacun} Let  $f \in  L^p_{loc}$ in a neighborhood of $x_0$ for a $p >1$, and assume that the $p$-exponent of $f$ is finite in  a left neighborhood of $p_0(f) $. The {  lacunarity exponent } of $f$ at $x_0$  is
 \BE  \label{deflac}   {\mathcal L}_f (x_0) = \frac{\partial }{\partial q} \left( h_f^{1/q} (x_0)\right)_{ q=q_0 (f)^+} . \EE 
 \ED
 
 { \bf Remarks:}
 
 \begin{itemize}
\item Even if the $p$-exponent is not defined at $p_0 (f)$, nonetheless, because of the concavity  of the mapping $q \rightarrow h^{1/q}_f (x_0)$ (see Proposition   \ref{propinterp} below), its right derivative is always well-defined, possibly as a limit.  
\item As expected, the lacunarity exponent of a cusp vanishes, whereas the lacunarity exponent of a lacunary comb coincides with  the accessibility exponent   of its support. 
\item The condition  $ {\mathcal L}_f (x_0) \neq 0$ does not mean that  the support of $f$ (or of $f-P$) has a positive accessibility exponent (think of the function $F^\al_{ \ome, \gamma} +g $ where $g$ is a $C^\infty$  but nowhere polynomial function).  
\item  The definition supplied by (\ref{deflac})  bears similarity with the definition of the oscillation exponent (see \cite{Oscillsing,Bergou}  and ref. therein) which is also defined through a derivative of a pointwise exponent;  but the variable with respect to which the derivative is computed  is the order of a fractional integration. The relationships  between these two exponents  will be investigated in a forthcoming paper \cite{PART1}.
\end{itemize}

  \section{\boldmath Properties of the $p$-exponent}
 
 \label{fracexp} 
 
 In signal and image processing, one often meets data that cannot be modeled by functions  $f \in L^1_{loc}$, see Figure \ref{fig:rwimage}.  
 It is therefore necessary to set the analysis in a wider functional setting, and therefore to extend the notion of $T^p_\al (x_0)$ regularity to the case $p <1$.  
 
 \subsection{\boldmath The case $p <1$}

 \label{subsecplone}

 The standard way to perform this extension  is to consider exponents in the setting of the real Hardy spaces $H^p$ (with $p <1$)  instead of $L^p$ spaces, see \cite{JaffCies,JaffToul}.     First, we  need to extend the definitions that we gave  to the range   $p \in (0, 1]$. The simplest way   is to start with the wavelet characterization of $L^p$ spaces, which we now recall. 
 
 We denote indifferently by $\chi_{j,k} $ or  $\chi_\la$ the characteristic function of the interval $\la\;  (= \la_{j,k})$ defined by (\ref{deflamb}). The { \bf wavelet square function }  of $f$ is 
\[{ \mathcal W}_f (x) =  \left(\sum_{(j,k)\in \ZZ^2} | \cli |^2 \chi_{j,k} (x)  \right)^{1/2}. \] 
Then, for $p >1$,  
\BE \label{aphp} f \in L^p (\RR ) \Longleftrightarrow  \int_\RR\left({\mathcal W}_f (x)\right)^p dx < \infty ,  \EE
see \cite{Mey90I}.  The quantity 
$\left(  \int \left({\mathcal W}_f (x)\right)^p dx \right)^{1/p} $   is
thus  equivalent to $\parallel f \parallel_p$. One  can then take the characterization supplied by (\ref{aphp}) when $p >1$ as a definition of the Hardy space $H^p$ (when $p \leq 1$); note that this definition yields  equivalent quantities when the (smooth enough) wavelet basis is changed, see \cite{Mey90I}.  This justifies the fact that we will
 often denote  by $L^p$ the space $H^p$, which will lead to no confusion;  indeed, when $p \leq 1$ this notation will refer to $H^p$, and,  when $p  > 1$ it  will refer to $L^p$.    

 Note that, if $p=1$, (\ref{aphp}) does not characterize the space $L^1$ but a strict subspace of $L^1$ (the real Hardy space $H^1$, which consists of  functions of $L^1$ whose Hilbert transform also belongs to $L^1$, see \cite{Mey90I}).

 Most results proved for the $L^p$ setting will extend without modification to the $H^p$ setting. In particular, $T^p_\al$ regularity can be extended to the case $p \leq 1$   and has the same wavelet characterization, see \cite{JaffTPu}. All definitions introduced previously therefore extend to this setting.  \\  


The definition of $T^p_\al (x_0)$ regularity given by (\ref{pexpa}) 
 is  a size estimate of an $L^p$ norm restricted to intervals $[x_0 -r, x_0 +r]$. Since  the elements of $H^p$  can be distributions,  the restriction of $f$ to an interval cannot be done directly (multiplying a distribution by a non-smooth function, such as a characteristic function, does not always make sense). This problem can be solved as follows:  If $I$ is an open interval,  one defines
$ \parallel f\parallel_{ H^p (I )} = \inf  \parallel g\parallel_p$, where the infimum is taken on the $g\in H^p$ such that $f = g $ on $I$. 
The  $T^p_\al$ condition  for $p \leq 1$ is then defined by: 
\[  f \in T^p_\al (x_0)  \quad \Longleftrightarrow \quad \parallel f\parallel_{ H^p ((x_0- r, x_0 +r) )} \leq C\; r^{\al +1/p},\] also when $p<1$.
We will show below  that  the $p$-exponent takes values in $[ -1/p, + \infty  ]$. 


\subsection{\boldmath When  can one use $p$-exponents?}

We already mentioned that, in order to use the H\"older exponent as a way to measure pointwise regularity, we need to check that the data are locally bounded, a condition which is  implied by the criterion $\Hmin >0$, which is therefore used as a practical prerequisite. Similarly, in order to use  a $p$-exponent based multifractal analysis, we need to check that the data locally belong to $L^p$ or $H^p$, a condition which can be  verified in practice  through the computation of the { \em wavelet scaling function}, which we now recall.

    The Sobolev space $L^{p,s}$ is defined by 
    \[ \forall s \in \RR, \; \forall p >0, \qquad  f \in L^{p,s} \quad \Longleftrightarrow \quad (Id -\Delta )^{s/2} f \in L^p , \]
    where the operator $(Id -\Delta )^{s/2}$ is the Fourier multiplier by $(1+ | \xi |^2)^{s/2} $, and we 
    recall our convention that $L^p$ denotes the space $H^p$ when $p \leq1$, so that Sobolev spaces are defined also for $p \leq 1$. 
    
    \BD 
    Let $f$ be a tempered distribution. The wavelet scaling  function of $f$ is defined by 
\BE \label{scalsob} \forall p >0, \qquad \eta_{f} (p) =p  \;\sup \{ s: f \in L^{p,s}  \} . \EE
\ED

Thus, $\forall p >0$: 

\begin{itemize}
\item If $\eta_{f} (p) >0$ then $f \in L^p_{loc}$.
\item If $\eta_{f} (p) < 0$ then $f \notin L^p_{loc}$.
\end{itemize}

 The wavelet characterization of Sobolev spaces implies that   the   wavelet scaling function  can be expressed as (cf. \cite{jmf2})
  \BE \label{defscalond}  \forall p >0, \hspace{6mm} 
  \eta_f (p) =   \displaystyle \liminf_{j \rightarrow + \infty} \;\; \frac{\log \left( 2^{-j} \displaystyle\sum_{ k }  | c_{j,k} |^p \right) }{\log (2^{-j})}. \EE 
  
This provides a practical criterion for determining if data locally belong to $L^p$, supplied by the condition $\eta_{f} (p) >0$. 
The following bounds for   $p_0 (f)$ follow:
\[ \sup \{ p : \eta_{f} (p) >0\}  \leq p_0 (f)  \leq \inf \{ p : \eta_{f} (p) <0 \}, \]
which (except in the very particular cases where $\eta_{f}$ vanishes identically on an interval) yields the exact value of $p_0 (f)$. 

\begin{figure}[tb]
\centering
\includegraphics[width=\linewidth]{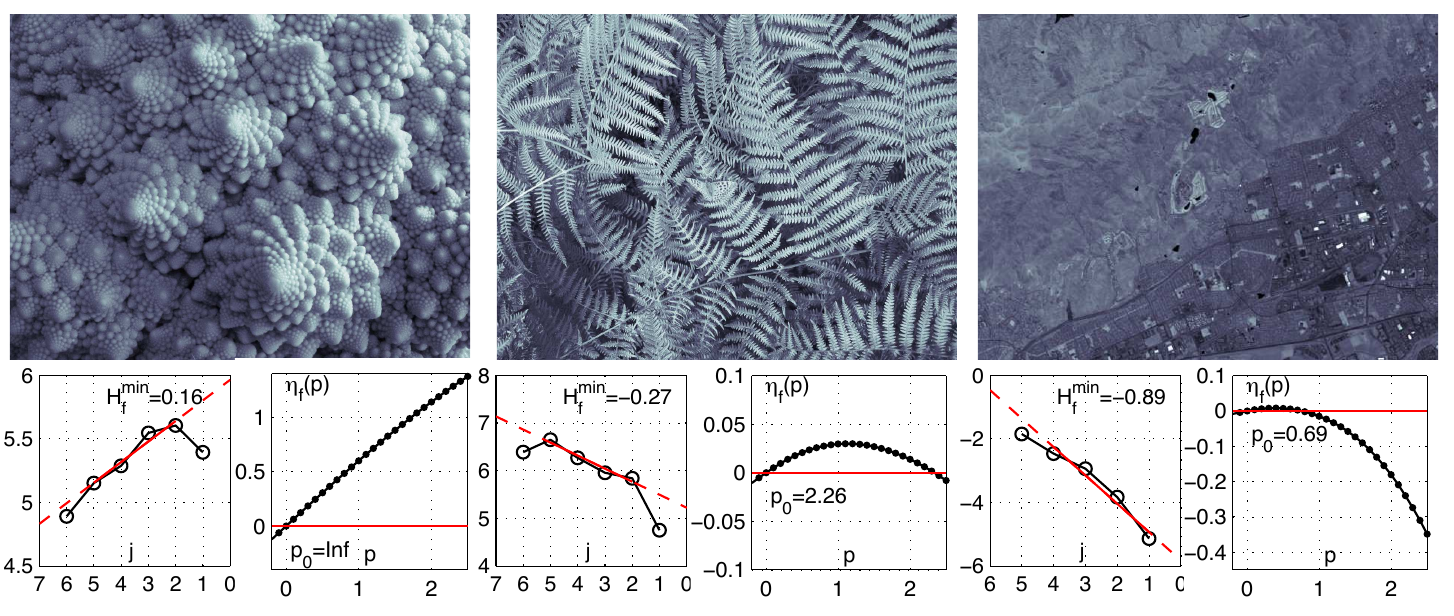}
\caption{\label{fig:rwimage}{\bf Real-world images} (top row) of Romanesco
  broccoli (left column), fern leaves (center column) and a patch of a hyperspectral image of the Moffett field, acquired by the AVIRIS instrument  (spectral band 90, right column).
Bottom row: corresponding numerical estimation of uniform H\"older exponents $H_f^{min}$, wavelet scaling functions $\eta_f(p)$ and critical Lebesgue indices $p_0$, respectively.}
\end{figure}

In applications, data with very different values of $p_0 (f)$ show up; therefore, in practice, the mathematical framework supplied by the whole range of $p$ is relevant.  
As an illustration, three examples of real-world images with positive and negative uniform H\"older exponents and with critical Lebesgue indices above and below $p_0=1$ are analyzed in Figure \ref{fig:rwimage}.

\setlength{\tabcolsep}{1pt}

\subsection{\boldmath Wavelet characterization  of $p$-exponents  }

In order to compute and prove properties of $p$-exponents we will need the exact wavelet characterization of $T^p_\al (x_0)$, see \cite{JaffTPu, JaffCies}. 
Let $\lambda$ be a dyadic interval; $3 \lambda$ will denote the interval  of same center and three times wider (it is the union of $\la $ and its two closest neighbors).  For $x_0 \in \RR^d$,  denote by $\lambda_j (x_0)$  the dyadic cube of width  $2^{-j} $  which contains $x_0$.  The { \bf local square functions}  at $x_0$ are  the sequences defined for $j \geq 0$ by  
\[{ \mathcal W}^j_{f , x_0} (x) =  \left(\sum_{\la \subset 3 \lambda_j (x_0)} | \ c_\la |^2 \chi_\la (x)  \right)^{1/2}. \] 
Recall that (cf. \cite{JaffTPu})
\BE \label{caractpal}  f\in T^p_\al (x_0) \quad \mbox{ if and only if } \exists C>0,\;\forall j\geq 0\qquad
  \left\| { \mathcal W}^j_{f , x_0}   \right\|_{p} \leq C\;2^{-(\al +1/p)j}.\EE

The following result  is required for the definition of the lacunarity exponent in (\ref{deflac}) to make sense,  and implies that Definition \ref{deflacun} also makes sense  when $p_0 (f) <1$. 

\BP \label{propinterp} Let $p,q \in (0, + \infty ]$, and  suppose that $f \in T^p_\al (x_0) \cap T^q_\be (x_0)$;  let $\theta \in [0,1]$. Then  $f \in T^r_\gamma (x_0) $, where 
\[  \frac{1}{r} =  \frac{\theta}{p} + \frac{1-\theta}{q} \quad \mbox{ and} \quad \gamma = \theta \al + (1-\theta ) \be. \]
It follows that the mapping $q\rightarrow h^{1/q}_f (x_0)$ is concave on its domain of definition. 
\EP

{ \bf Proof:}  When $p, q < \infty$, the result is  a consequence of (\ref{caractpal}).  
H\"older's inequality implies that 
\[  \left\| { \mathcal W}^j_{f , x_0}   \right\|_{r} \leq   \left\| { \mathcal W}^j_{f , x_0}   \right\|^{\theta/p}_{p}   \left\| { \mathcal W}^j_{f , x_0}   \right\|^{(1-\theta)/q}_{q}  .\] 
We thus obtain the result  for $p, q  < \infty$. The case when $p$ or $q = + \infty$  does not follow, because  there exists no exact wavelet characterization of $C^\al (x_0)= T^\infty_\al (x_0)$; however, when $p, q >1$,  one can use the initial  definition of $  T^p_\al (x_0)$ and $C^\al (x_0)$  through local $L^p$ and $L^\infty$ norms  and the result also follows from H\"older's inequality; 
hence  Proposition \ref{propinterp} holds. \\

If $f\in H^p$, then $\parallel { \mathcal W}_f  \parallel_p \leq C$. Since $ { \mathcal W}^j_f  \leq  { \mathcal W}_f $, it follows that  $\parallel { \mathcal W}^j_f  \parallel_p \leq C$, so that (\ref{caractpal}) holds with $\al =-1/p$. Thus $p$-exponents are always larger than $-1/p$ (which extends to the range $p <1$ the result already mentioned for $p \geq 1$).   Note that this bound is compatible with the existence of  singularities of arbitrary large negative order (by picking $p$ close to 0). The example of cusps will now show that the $p$-exponent  can indeed take values down to $-1/p$.

\setlength{\tabcolsep}{1pt}
\begin{figure}[tb]
\includegraphics[width=\linewidth]{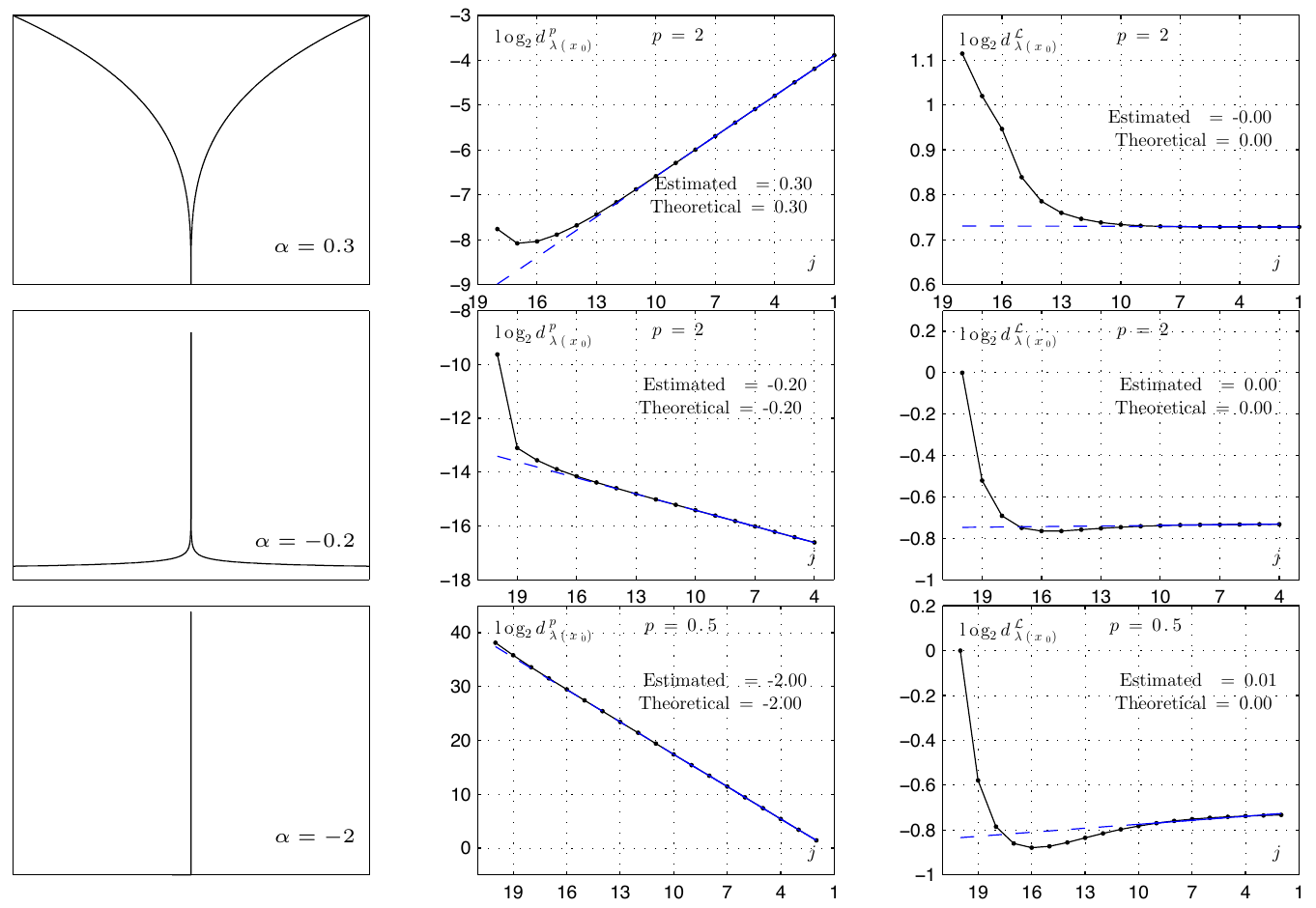}
\caption{\label{fig:cusp}{\bf Cusps} with exponents $\alpha=\{+0.3,\;-0.2,\;-2\}$  (from top to bottom row, respectively): functions (left column) and estimation of $p$-exponents and lacunarity exponents (center and right column, respectively). The critical Lebesgue indices are given by $p_0=\{+\infty,\;5,\;0.5\}$, respectively.}
\end{figure}

\subsection{\boldmath Computation of $p$-exponents for cusps }

Typical examples of distributions for which the $p$-exponent is constant (see Proposition \ref{pregcusps} below) and equal to a given value $\al < -1$  are supplied by the cusps ${ \mathcal C}_{\al}$,   whose definition  can be extended to the range  $\al \leq -1$ as follows: First, note that cusps cannot  be defined directly for $\al \leq -1$  by (\ref{cusp}) because they do not belong to $L^1_{loc}$  so that  they would be ill-defined even in the setting of distributions (their integral against a $C^\infty$ compactly supported function $\varphi $ may diverge). Instead, we use the fact that, 
if $ \al >1$, then $  { \mathcal C}_{\al}''= \al (\al -1) { \mathcal C}_{\al-2}$,  which indicates a way to define by recursion the cusps ${ \mathcal C}_{\al} $, when $\al <-1$ and $\al \notin \ZZ$, as follows: 
\[ \mbox{ if } \quad \al <0, \qquad  { \mathcal C}_{\al} =   \frac{1}{(\al +1)(\al +2)}{ \mathcal C}_{\al+2}'' , \]
where the derivative is taken in the sense of distributions. 
The ${ \mathcal C}_{\al}$ are thus defined as distributions when $\al$ is not a negative integer.  It can also be done  when $\al $ is a negative  integer, using the following definition for $\al =0$ and $ -1$:
\[  { \mathcal C}_{0} = \log (|x| )\qquad \mbox{ and}  \qquad { \mathcal C}_{-1} = { \mathcal C}_{0}' = P.V. \left(\frac{1}{x}\right) , \] 
where P.V. stands for ``principal value''. 
%

\BP \label{pregcusps}
 If $\al \geq 0$,  the cusp ${ \mathcal C}_{\al}$ belongs to $L^\infty_{loc} $   and  its $p$-exponent is $\al$. 
  If $\al  <  0$,  the cusp ${ \mathcal C}_{\al}$ belongs to $L^p_{loc} $ for $ p < -1/\al$   and  its $p$-exponent is $\al$. 
  \EP
  

  { \bf Proof of Proposition \ref{pregcusps}:} The case  $\al \geq 0$  and $p \geq 1$ has already been considered in \cite{PART1,JaffMel}.  In this case, the computation of the $p$-exponent is straightforward.
  Note that, when $\al \in (-1, 0)$ and $p \geq 1$ the computations are similar.  We thus focus on the distribution case, i.e.  when $p < 1$. The global and pointwise regularity will be determined through an estimation  of the wavelet coefficients of the cusp.  We use a smooth enough, compactly supported wavelet basis and we denote by $c_{j,k}$ the wavelet coefficients of the cusp
  \[ c_{j,k} = 2^j \langle \psi_{j,k} | { \mathcal C}_{\al} \rangle . \]
  The selfsimilarity of the cusp implies that 
  \BE \label{cuspwc1}\forall j,k \qquad c_{j,k}  = 2^{-\al j} c_{0,k} ; \EE
  additionally, as soon as $k$ is large enough so that  the support of $\psi (x-k)$ does not intersect the  origin, the cusp is $C^\infty$ in the support of $\psi (x-k)$  and coincides with the function $| x|^\al$. An integration by parts then yields that, for any $N$ smaller than the global regularity of the wavelet,
  \[ c_{0,k}  = (-1)^N \int {\psi^{(-N)}} (x-k) \; \al (\al-1) \cdots (\al-N) |x|^{\al-N}  dx ,\] 
  so that the sequence $c_{0,k}$ satisfies 
  \BE \label{cuspwc2} |c_{0,k} | \leq  \frac{C_N}{(1+|k|)^N}\EE
  where $N$ can be picked arbitrarily large. 
  The estimation of the $L^p$ norm of the wavelet square function follows easily from    (\ref{cuspwc1}) and (\ref{cuspwc2}), and so does the lower bound for the $p$-exponent. The upper bound is obtained by noticing that one of the $c_{0,k}$ does not vanish (otherwise, all  $c_{j,k}$ would vanish, and the cusp would be a smooth function at the origin).   Therefore, there exists at least one $k_0$ such that $\forall j $, $c_{j,k_0} = C 2^{-\al j}$, and the wavelet characterization of $T^p_\al$ regularity then yields that $h^p (x_0) \leq \al$. \\
  
Three examples of cusps and numerical estimates of their $p$-exponents and lacunarity exponents are plotted in Figure \ref{fig:cusp}.


%
   

 \subsection{Wavelet characterization and thin chirps}
 
 \label{thinch}
 
 In practice, we will derive   $T^p_\al$ regularity from simpler quantities than the local square functions. 
 The { \bf\boldmath $p$-leaders} of $f$ are defined by local $l^p$ norms of wavelet coefficients as follows:
 \begin{equation} \label{qleaders} d_\lambda^p = \left( \sum_{ \lambda' \subset 3 \lambda} | c_{\lambda'} |^p 2^{-(j'-j)} \right)^{1/p}  \end{equation}
(they are finite if $f \in L^p_{loc} (\RR^d)$, see \cite{JaffMel}). Note that, if $p = + \infty$, the corresponding  quantity $d_\lambda^\infty$  
is usually denoted by $d_\lambda$ and simply called the 
  { \bf wavelet leaders}; we have 
 \BE \label{defwl} d_\lambda := d_\lambda^\infty  = \sup_{\lambda ' \subset 3 \lambda} | c_{ \lambda '} |  . \EE

The  notion of   $T^p_\alpha $ regularity can be related to $p$-leader coefficients (see \cite{JaffToul, Bergou, JaffMel}):

 \begin{equation} \label{carachqf} \mbox{ If $\eta_f (p) >0$, then  } \quad h^{p}_f (x_0)=  \liminf_{j \rightarrow + \infty}\;  \; 
  \frac{ \log 
 \left( d^p_{\lambda_j (x_0)} \right)  }{\log (2^{-j})}.
\end{equation}


\begin{figure}[tb]
\includegraphics[width=\linewidth]{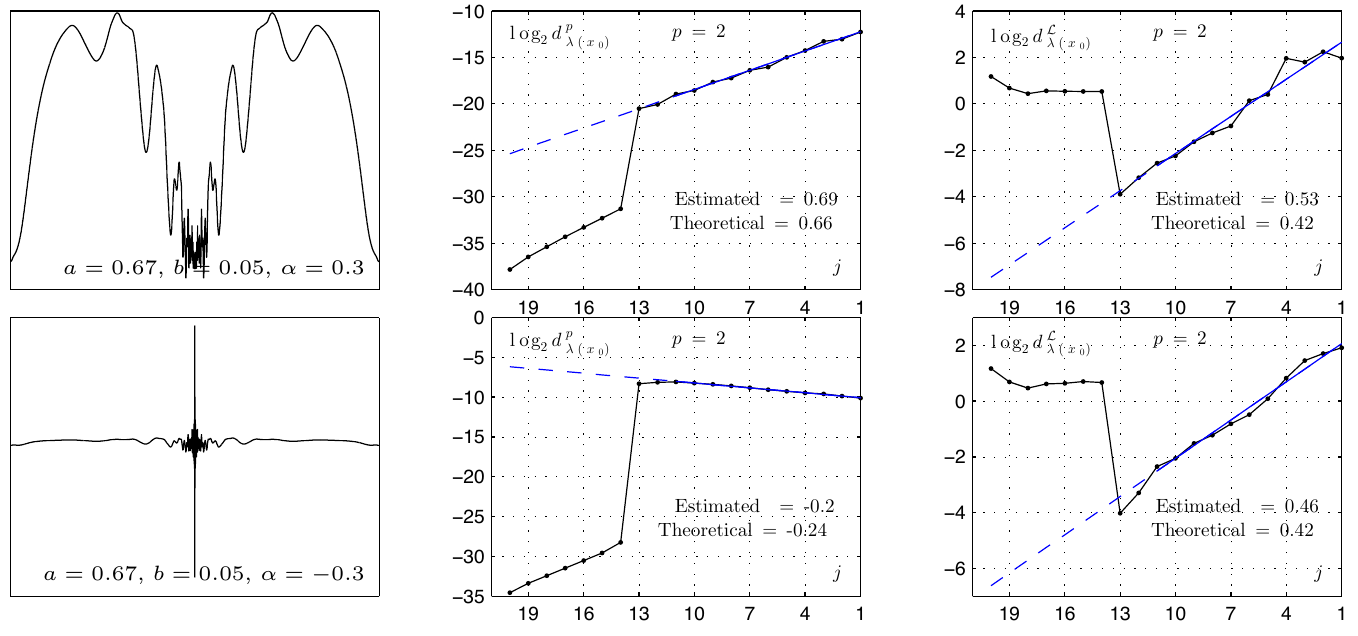}
\caption{\label{fig:Tchirp}{\bf Thin chirps} with $p_0=\infty$ (top row) and $p_0=3.2$ (bottom row): functions (left column) and estimation of $p$-exponents and lacunarity exponents (center and right column, respectively).}
\end{figure}

Our purpose in this section is to   introduce new  ``archetypical'' pointwise singularities which will yield examples where the $p$-exponent and the lacunarity exponent  can take arbitrary values.   Because of (\ref{carachqf}),   it is easier to work with examples that are defined directly by their wavelet coefficients on a smooth wavelet basis. We therefore develop new examples rather than extending the lacunary combs of Section \ref{sec:lacunary}. 

 \BD \label{defithinchir} Let $a, b \in (0, 1)$ satisfying $0 < b<1-a$, and let $\al \in \RR$. The {thin chirp}  ${ \mathcal T}_{a,b, \al}$ is defined by its wavelet series 
 \[ { \mathcal T}_{a,b, \al} = \sum_{ j \geq 0} \; \sum_{k \in \ZZ} c_{j,k} \; \psi_{j,k}, \]
 where
\[  \begin{array}{rll} c_{j,k} & = 2^{-\al j}   &  \mbox{ if } \quad k \in [  2^{(1-a)j},  2^{(1-a)j} +  2^{bj} ] \\ 
 & = 0  &   \mbox{ otherwise.  } \end{array}    \]
\ED

The following results are straightforward, using the wavelet characterization of $L^p$ and $T^p_\al$ regularity.  

\BP
The {thin chirp}  ${ \mathcal T}_{a,b, \al}$  is bounded if and only if $\al  > 0$. 
\[  \mbox{ If} \quad \al \leq 0, \qquad  p_0 ({ \mathcal T}_{a,b, \al}) = \frac{1-b}{-\al} . \]

The $p$-exponent of ${  \mathcal  T}_{a,b, \al}$ at the origin  is 
\[  h^p_{{\mathcal  T}_{a,b, \al}}(0) =  \frac{1-a-b}{a} q  + \frac{\al}{a}  .    \]
\EP


Note that, if the wavelets are compactly supported, then for $j$ large enough the pack  of $2^{bj}$ successive wavelets with non-vanishing coefficients covers an interval of length $2^{-j} 2^{bj}$ at  a distance $2^{-aj}$ from the origin, so that the accessibility exponent of the support of  ${  \mathcal  T}_{a,b, \al}$ is $(1-a-b)/a$: Thus, it coincides with  the lacunarity exponent of ${ \mathcal T}_{a,b, \al}$  as expected.  \\

Illustrations of thin chirps and  the numerical estimation of their $p$-exponents and lacunarity exponents are provided in Figure \ref{fig:Tchirp}.

\subsection{\boldmath $p$-exponent analysis of measures}

Several types of measures (such as multiplicative cascades)  played a central role in the development of multifractal analysis. Since measures (usually) are not $L^1$ functions, their $p$-exponent for $p \geq 1$ is not  defined. Therefore,  it is natural to wonder if it can be the case when $p <1$.  This is one of the purposes of Proposition \ref{propregmes}, which yields 
sufficient conditions  under which  a measure $\mu$ satisfies $\eta_\mu (p) >0$ for $p <1$, which will imply that its $p$-exponent multifractal analysis can be performed. 
An important by-product of using $p$-exponents for $p \leq1$ is that it offers a common setting to treat pointwise regularity of measures and functions.

 Recall that $\overline{\dim_B}  (A)$ denotes the upper box dimension of the set $A$.

\BP \label{propregmes} Let $\mu$ be a measure; then its wavelet scaling function  satisfies  
 $\eta_\mu (1) \geq 0$. Furthermore,   if $\mu$ does not have a density which is an $L^1$ function, then $  \eta_\mu (1)=  0 . $

Additionally, if $\mu$ is a singular measure whose support $ supp (\mu )$  satisfies 
\BE \label{hypsupp}  \delta_\mu := \overline{\dim_B}  ( supp (\mu )) < 1, \EE
 then 
\BE \label{minregglob}  \forall p <1, \qquad  \eta_\mu (p) \geq {(1-\delta_\mu)}{(1-p)} , \EE
and
\BE \label{miaxregglob}  \forall p >1, \qquad  \eta_\mu (p) \leq {(1-\delta_\mu)}{(1-p)} . \EE
\EP

{ \bf Remarks:}
\begin{itemize}
\item  (\ref{minregglob}) expresses the fact that, if $\mu$ has a small support, then its Sobolev regularity is increased for $p <1$. This is somehow counterintuitive, since one expects a measure to become more singular when the size of its support shrinks; on the other hand (\ref{miaxregglob}) expresses that this is actually the case when  $p >1$. 
\item Condition $\delta_\mu <1$ is satisfied  if $\mu$ is supported by a Cantor-like set, or by a selfsimilar set satisfying Hutchinson's open set condition. 
\item (\ref{minregglob})  has an important consequence for the multifractal analysis of measures: Indeed, if $\delta_\mu < 1$,  then $\eta_\mu (p)  >0$ for $
p <1$, so that the classical mathematical results concerning the multifractal analysis based on the $p$-exponent apply, see Section \ref{seclacu}. 
\item A slightly different  problem was addressed by H. Triebel: In \cite{Trie},  he determined under which conditions the scaling functions  commonly used in the multifractal analysis of probability measures (see \eqref{defscalondnew} below) can be recovered through Besov or Triebel-Lizorkin norms (or semi-norms).
\end{itemize}

{ \bf Proof of Proposition \ref{propregmes}:} If $\mu$ is a measure, then for any continuous bounded function $f$
\BE \label{boumes} | \langle \mu | f \rangle | \leq C \parallel f \parallel_\infty. \EE
We pick 
\[ f = \sum_k \ep_{j,k} \psi_{j,k} , \qquad \mbox{ where} \quad \ep_{j,k}  = \pm 1, \]
so that $f$ is continuous and satisfies $\parallel f \parallel_\infty\leq C$, where $C$ depends only on the wavelet (but not on the choice of the $\ep_{j,k} $). Denoting by $c_{j,k}$ the wavelet coefficients of $\mu$, we have 
\[ \langle \mu | f \rangle  = \sum_k \ep_{j,k} \int \psi_{j,k}  d \mu = 2^{-j} \sum_k \ep_{j,k}c_{j,k} . \]
Picking $\ep_{j,k}= \mbox{sgn} (c_{j,k})  $ it follows from (\ref{boumes})  that 
\BE \label{bes1} 2^{-j} \sum_k | c_{j,k} | \leq C, \EE
or, in other words,  $\mu$ belongs to the Besov space $ B^{0,\infty}_1$, which implies that   $\eta_\mu (1) \geq 0 $, see \cite{Jaffard2004, Mey90I}.   

On  other hand, if $\mu \notin L^1$, then using the interpretation of the scaling function in terms of Sobolev spaces given by (\ref{scalsob}),  we obtain that $\eta_\mu (1) \leq 0$. Hence the first part of the proposition holds.

We now prove (\ref{minregglob}).  We assume that the used wavelet is compactly supported, and that its support is included in the interval 
$[ -2^l, 2^l]$ for an $l >0$ (we pick the smallest $l$ such that this is  possible).  Let $\delta > \overline{\dim_B}  ( supp (\mu )) $;   for $j$ large enough, $supp (\mu)$ is included in at most     $2^{[\delta j]    }$  intervals  of length $2^{-j}$.  It follows that, at scale $j$,  there exist at most   $2^{[\delta j]    } \cdot 2 \cdot 2^l$ wavelets $(\psi_{j,k} )_{k\in \ZZ}$  whose support intersects the support of $\mu$. Thus for  $j$ large enough, there are at most $C 2^{\delta j}$ wavelet coefficients that do not vanish.

 Let $p \in (0,1)$, $q=1/p$ and $r$ be the conjugate exponent of $q$, i.e. such that $1/q + 1/r =1$. Using H\"older's inequality,  
\[ \sum_k | c_{j,k}|^p \leq \left(  \sum_k | c_{j,k}|^{pq}\right)^{1/q} \left(  \sum_k 1^{r}\right)^{1/r}, \]
where the sums are over at most $C 2^{\delta j}$ terms; thus 
\[ \sum_k | c_{j,k}|^p \leq \left(  \sum_k | c_{j,k}|\right)^p  \; C \, 2^{\delta j/r}. \]
Using (\ref{bes1}), we obtain that 
\[  2^{-j} \sum_k | c_{j,k}|^p \leq  C 2^{-(1-\delta)j/r},\] 
so that $\eta_\mu (p)  \geq (1-\delta) ( 1-p) $. Since this is true  $\forall \delta > \delta_\mu$, (\ref{minregglob}) follows. 

We now prove (\ref{miaxregglob}).  Let $p \geq 1$  and let $q$ be the conjugate exponent. Using H\"older's inequality, 
\[ \sum_{k} |c_{j,k} |  \leq \left(\sum_{k} |c_{j,k} |^p \right)^{1/p} \left(\sum_{k} 1^q\right)^{1/q} . \]
Let again $\delta > \delta_\mu$; using the fact that the sums bear on at most $2^{\delta j}$ terms,  and that the left-hand side is larger than $C2^j$, we obtain that 
\[   \left(\sum_{k} |c_{j,k} |^p \right)^{1/p}  \geq C\; 2^{ j} 2^{-\delta j/q},\]
which can be rewritten
\[   2^{- j} \sum_{k} |c_{j,k} |^p  \geq C\; 2^{ -j} 2^{ pj}  2^{-\delta jp/q},\]
  so that 
  $ \eta_\mu (p) \leq (1-p) (1-\delta)$; since this is true  $\forall \delta > \delta_\mu$, (\ref{miaxregglob}) follows,  
  and  Proposition \ref{propregmes} is completely proved. \\

Since $p=1$ is a borderline case for the use of the $1$-exponent one may expect that picking $p <1$ would yield $\eta_\mu (p) > 0$ (in which case one would be on  the safe side in order to recover mathematical results concerning the $p$-spectrum, see \cite{JaffCies, Abel}). 
However, this is not  the case, since there exist even continuous functions $f$ that satisfy $\forall p >0$,  $\eta_f (p)= 0$. An example is supplied by  
\[ f = \sum_{j \geq 0} \sum_{k \in \ZZ} \;  \frac{1}{j^2}\,  \psi_{j,k} . \]


%


\section{ Multifractal analysis of lacunary wavelet series }

\label{seclacu}

 Multifractal analysis is motivated by the observation that many mathematical models have an extremely erratic pointwise regularity  exponent  which jumps everywhere; this is the case e.g. of  multiplicative cascades or of L\'evy processes, whose exponents $h$  satisfy  that
\BE \label{oscil}  \mbox{ a.s.} \quad \forall x_0, \qquad \limsup_{x_\rightarrow x_0} h(x)- \liminf_{x_\rightarrow x_0} h(x) \EE
is bounded from below by a fixed positive quantity (we will see that this is also the case for lacunary wavelet series).  This clearly excludes the possibility of any robust direct estimations of $h$.  The driving idea of multifractal analysis is that one should rather focus on alternative quantities that
\begin{itemize}
\item  are numerically computable on real life data in a stable way,
\item  yield  information on the erratic behavior of the pointwise exponent.
\end{itemize}
Furthermore, for standard random models (such as the ones mentioned above)  we require these quantities not to be random (i.e. not to depend on the sample path which is observed) but to depend  on the characteristic parameters of the model only. 
The relationship between the { \bf multifractal spectrum} and scaling functions (initially pointed out by U. Frisch and G. Parisi in \cite{ParFri85}; see \eqref{formult1} below) satisfies these requirements.

We now recall the notion of multifractal spectrum. We denote by $\dim (A)$ the Hausdorff dimension of the set $A$. 

\BD
Let $h (x)$ denote a pointwise exponent. 
The multifractal  spectrum $d (H)$ associated with this pointwise exponent  is 
\[ d(H) =  \dim \{ x: \;\; h (x)=H\} . \] 
\ED

In the case of the $p$-exponent, the sets of points with  a given $p$-exponent will be denoted by $ F^p_f (H)$:
\BE \label{givhpex}  F^p_f (H) = \{ x_0: \; h^p_f(x_0) = H\} , \EE 
and the corresponding  multifractal spectrum (referred to as the $p$-spectrum) is denoted by $d^p (H)$;  in the case of  the lacunarity exponent, we denote it by $d^\mathcal{L} (L) $. 

\subsection{Derivation of the multifractal formalism}

\label{derivfp}

We now recall how  $d (H)$  is expected to be  recovered from global quantities effectively computable on real-life signals (following the seminal work of G. Parisi and U. Frisch \cite{ParFri85} and its wavelet leader reinterpetation \cite{Jaffard2004}). 
A key assumption is that this exponent  can be  derived from nonnegative quantities (which we denote either by $e_{j,k}$ or  $e_\la$), which are defined on the set of dyadic intervals, by a log-log plot regression: 
 \begin{equation} \label{caracexppoi} h (x_0)=  \liminf_{j \rightarrow + \infty}\;  \; 
  \frac{ \log 
 \left( e_{\lambda_j (x_0)} \right)  }{\log (2^{-j})}.
\end{equation}
It is for instance the case of the $p$-exponent, as stated in (\ref{caractpal}) or (\ref{carachqf}), for which the quantities $e_{\lambda}$ are given by the $p$-leaders $d^p_{\lambda}$.

In the case of the lacunarity exponent,  quantities $e_{\lambda}$ can be derived as follows: Let $\Delta q >0$ small enough be given. 
 If $f$ has a $1/q$-exponent $H$ and  a lacunarity exponent  $L$ at $x_0$ then its $1/q$-leaders satisfy
\[ d^{1/q}_{j }  (x_0)  \sim 2^{-Hj}  ,\]
 and its $1/(q + \Delta q)$-leaders satisfy
 \[ d^{1/(q + \Delta q)}_{j}(x_0)    \sim 2^{-(H+ \Delta q L)j} ;\]
 we   can eliminate $H$ from these two quantities by considering the { \bf ${\mathcal{L}} $-leaders:}
\[ d^{\mathcal{L}}_{\lambda} : = \left( \frac{d^{1/(q + \Delta q)}_{j}}{d^{1/q}_{j }  }  \right)^{1/\Delta q}  \sim 2^{- Lj} .\] 
(this argument follows a similar one developed in \cite[Ch. 4.3]{Bergou} for the derivation of a  multifractal analysis associated with the oscillation exponent).  \\

The multifractal spectrum will be derived from  the following quantities, referred to as the { \em structure functions}, which are similar to the ones that come up in the characterization of the wavelet scaling function in (\ref{defscalond}): 
 \[ S_{f} (r,j) = \left( 2^{-j} \displaystyle\sum_{ k }  | e_{j,k} |^r \right).  \] 
The  scaling function associated with the collection of $(e_\la)$ is 
  \BE \label{defscalondnew}  \forall r \in \RR , \hspace{6mm} 
  \zeta_f (r) =   \displaystyle \liminf_{j \rightarrow + \infty} \;\; \frac{\log \left( S_{f} (r,j)  \right) }{\log (2^{-j})}. \EE 
Let us now sketch the heuristic derivation of the multifractal formalism; (\ref{defscalondnew})  means 
 that, for  large $j$,
 \[  S_{f} (r,j)   \sim 2^{- {\zeta} (r)j}.  \]
Let us  estimate the contribution to $ S_{f} (r,j)  $ of the dyadic intervals  $  \la$ that cover the points of $E_H$.
By definition of $E_H$, they satisfy
$ e_\la  \sim 2^{-Hj}; $
 by definition of $d(H)$,   since we use cubes of the same width $2^{-j}$  to cover $E_H$, we need about $2^{d (H) j}$ such cubes; therefore the corresponding
contribution   is  of the order of magnitude of 
\[ 2^{-j}2^{d (H) j}2^{-Hrj} = 2^{-(1-d(H) +Hr)j}.\]
 When $j \rightarrow +
\infty$, the dominant contribution comes from the smallest exponent, so that 
\begin{equation} \label{etamuleg} {\zeta} (r) = 
\inf_H (1- d(H)+ Hr).\end{equation}

By construction, the  scaling function 
 ${\zeta}(r)$  is  a concave function on $\RR$, see \cite{ParFri85,Jaffard2004,riedi03} which 
  is in agreement with the fact that the right-hand side of (\ref{etamuleg}) necessarily is a
concave function (as an infimum of a family of linear functions) no matter whether $d (H)$ is concave or not.
If $d(H)$  also is   a concave function, then the Legendre transform in (\ref{etamuleg})  can be inverted
(as a consequence of  the duality of convex functions),
which justifies the following assertion.

\BD A nonnegative sequence $(e_\la)$, defined on the dyadic intervals, follows the   multifractal formalism   if the associated  multifractal spectrum $d(H)$  satisfies 
\begin{equation} \label{formult1}
d(H) = \inf_{r \in \RR} (1- {\zeta} (r) + Hr).
\end{equation} 
\ED

The derivation  given above is not a mathematical proof, and the determination of the  range of validity  of
(\ref{formult1}) (and of its variants) is one of  the main mathematical problems concerning  multifractal analysis.  
If it  does not  hold in complete generality,  
    the multifractal formalism  nevertheless yields an upper bound of  the spectrum of singularities, see \cite{ParFri85,Jaffard2004,riedi03}: As soon as (\ref{caracexppoi}) holds, 
\[ d(H) \leq \inf_{r \in \RR} (1- {\zeta} (r) + Hr). \] 

In  applications,  multifractal analysis is often used  only as a classification  tool in order to discriminate
between several types of signals; then, one is not directly concerned with the validity of (\ref{formult1})
but  only with a precise  computation of the new { \em multifractal parameters} supplied by the scaling function, or equivalently its Legendre transform.  Note that  studies of multifractality for the $p$-exponent have been performed  by A. Fraysse  who proved genericity results of multifractality for functions in Besov or Sobolev spaces in \cite{Fraysse}.

\subsection{Description of the model and global regularity}

In this section, we extend to possibly negative exponents the  model of   lacunary wavelet series  introduced in \cite{Jaf6}. 
We assume that 
$
\psi$ is   a wavelet  in the Schwartz class (see however the remark after Theorem \ref{theospec},  which gives sufficient conditions of validity of the results of this section  when wavelets of limited regularity are used).  Lacunary wavelet series  depend on a { \bf lacunarity parameter}  $\eta \in (0,1)$ and a { \bf regularity  parameter}  $\al \in \RR$. 
At each scale $j\geq 0$, the process $X_{\al, \eta}$ has  exactly $[2^{\eta j} ] $ nonvanishing  wavelet coefficients on each interval $[l, l+1) $ ($l \in \ZZ$),  their common size  is $2^{-\al j}$, and their locations are picked at random: In each interval $[l, l+1)$ ($l \in \ZZ$), all drawings of  $[2^{\eta j} ] $  among the $2^j$ possibilities $\frac{k}{2^j} \in [l, l+1)$ have the same probability. Such a series is called a { \bf  lacunary wavelet series} of parameters $(\al, \eta)$. 
 Note that, since $\al$ can be arbitrarily negative,  $X_{\al, \eta}$ can actually be a random distribution of arbitrary large order. 
By construction
\[ H^{min}_{X_{\al, \eta}} = \al , \]
and, more precisely, the sample paths of  $X_{\al, \eta}$ are locally bounded if and only if $\al >0$. 
The case considered in \cite{Jaf6} dealt with $\al >0$, and was  restricted to the computation of H\"older exponents.  Considering $p$-exponents allows us   to extend the model to negative values of $\al$, and also to see how the global sparsity of the wavelet expansion (most wavelet coefficients vanish) is related with the pointwise lacunarity of the sample paths. Note that extensions of this model in different directions have been worked out in \cite{AJ02, Durand}

Since we are interested in local properties of the process $X$, we restrict our analysis to the interval $[0,1)$ (the results proved in the following clearly do not depend  on the particular interval which is picked); we can  therefore  assume that $k \in \{ 0, \cdots 2^j -1\}$.

We first determine how $\al$ and $\eta$ are related with the global regularity of the sample paths. 
The characterization (\ref{defscalond}) implies that the wavelet scaling function  is given by 
\BE \label{etala}  \forall p >0, \qquad \eta_{X_{\al, \eta}}  (p) =   \al p - \eta +1. 
\EE
It follows that 
\[  p_0 := p_0 (X_{\al, \eta}) = 
\begin{cases}
\frac{\eta-1}{\al} & \mbox{if } \alpha<0\\
+\infty & \mbox{if } \alpha>0.
\end{cases}
\] 
Note that $p_0  $ always exists and is strictly positive, even if $\al$ takes arbitrarily large negative values.
 We recover the fact that $p$-exponents allow us to deal with  singularities of arbitrarily large negative order. We will see  that this is a particular occurrence of a general result, see Proposition \ref{propspars}; the key property here is the sparsity of the wavelet series.

 \subsection{ Estimation of the \boldmath$p$-leaders  of $X_{\al, \eta}$} 
 
  An important step in the determination of the $p$-exponent  of sample paths of $X_{\al, \eta}$ at every point  is the estimation of their $p$-leaders. 
 We now  assume that $p < p_0  $, so that the sample paths of $X_{\al, \eta}$ locally belong to $L^p$ and the $p$-exponent of $X_{\al, \eta}$ is well-defined everywhere. 
Recall that the {  $p$-leaders}  are defined  by 

\BE \label{eqllam}  l_\la = \left( \sum_{ \la ' \subset 3\la}  | c_{\la'} |^p 2^{-(j'-j)}\right)^{1/p} . \EE

The derivation of  the $p$-exponent of $X_{\al, \eta}$ everywhere will be deduced from the  
estimation of  the  size of the $p$-leaders of $X_{\al, \eta}$.  A key result  is supplied by the following proposition, which states that   the size of the $p$-leaders of a lacunary wavelet series  is correctly estimated by the size of the first nonvanishing wavelet coefficient of smaller scale  that is met in the set $\{ \la ' : \;\;  \la' \subset 3 \la\} $. 

\BP \label{sizepleadlws} Let $\al \in \RR$, $\eta \in (0,1)$ and  let $X_{\al, \eta}$ be a  lacunary wavelet series of parameters $(\al, \eta)$;  
for each dyadic interval  $\la$ (of width $2^{-j}$), we define $j' $ ( $= j' (  \la) $)  as  the smallest random integer such that
\[ \exists \la' \subset 3 \la \;\; \mbox{ such that } | \la' | = 2^{-j'} \;\; \mbox{ and} \quad c_{\la '} \neq 0 . \]
Then, a.s. $\exists J$, $\exists C, C' >0 $ such that $\forall j \geq J$, $\forall \la $ of scale $j$
\[  \qquad  C   2^{-\al j'} 2^{-(j'-j)/p}  \leq l_\la  \leq C'    2^{-\al j'}  2^{-(j'-j)/p} j^{2/p}
\] 
\EP

 { \bf Proof:} This result will be  implied by the exponential decay rate $2^{-(j'-j)}$ that appears in the definition of $p$-leaders together with the lacunarity of the construction; we will show that  exceptional situations where this would  not be true  (as a consequence  of  local accumulations of nonvanishing coefficients) have a small probability and ultimately  will be  excluded by a Borel-Cantelli type argument.  We now make this argument precise. 
For that purpose, we will need to show that the sparsity of wavelet coefficients is uniform, which will be expressed by a uniform  estimate on the maximal  number of nonvanishing coefficients $c_{\la'}$ that can be found for $\la' $ (at  a given scale $j'$) included in a given interval $3  \la$. Such an estimate can be derived  by interpreting the choice of the nonvanishing wavelet coefficients in the construction of the model as a coarsening (on the dyadic grid) of  an { \bf empirical process}.  Let us now recall this notion, and the standard estimate on the increments of the empirical process that we will need. 

Let $N_j = [ 2^{\eta j} ] $ denote the number of nonvanishing wavelet coefficients at scale $j$.  We can consider that the corresponding dyadic intervals $\la $ have been obtained first by picking at random $N_j$ points in the interval $[0,1]$ (these points  are now $N_j$ independent
uniformly distributed random variables on $[0,1]$), and then by associating to each point the unique dyadic interval of scale $j$ to which it belongs.  Let $\PNT$  be  the process starting from 0 at $t=0$, which is piecewise constant and   which jumps by 1 at each  random point thus determined. The family  of processes 
\BE \ANT= \sqrt{N_j} \left(\frac{\PNT}{N_j} -t\right) \EE
is called an { \bf empirical process}  on $[0,1]$. 
The size of the increments of the empirical process  on  a given interval yields   information on the number of random points picked in this interval. If it is of length $l$, then the expected number of points is  $l [ 2^{\eta j} ]  $, and the deviation from this average can be  uniformly bounded   using the following result  of W. 
Stute        which is a particular case of Lemma 2.4 of  \cite{Stute}. 

\BL \label{l4} There exist two positive constants $C'_1$ and $C'_2$
  such that, if  $0 < l < 1/8$, $N_j l \geq 1$ and $8 \leq A
\leq C'_1 \sqrt{N_j l}$, 
\[ \PP \left( \sup_{|t-s| \leq l} | \ANT -\alpha^{j}_s | >
A\sqrt{l} \right)
\leq \frac{C'_2}{l} e^{ -A^2/64}. \]
\EL

Rewritten in terms of $\PNT$, this means that 
\BE \label{11bis} \PP  \left(     \sup_{|t-s| \leq l}  | \PNT -\PNS -N_j (t-s) |  > A  \sqrt{N_j l}  \right)   \leq    \frac{C'_2}{l} e^{ -A^2/64}. \EE
Recall that the assumption $\lambda'\subset 3\lambda$ implies that  $3\cdot 2^{-j}\geq 2^{-j'} $.
We will apply Lemma \ref{l4} differently for small values of $j' $ where the expected number of  nonvanishing coefficients $c_{\la'}$ that can be found for $\la' $ (at  a given scale $j'$) included in a given interval $ \la$ is very small, and the case of  large $j'$  where this number increases geometrically.  \\ 

We first assume that 
\BE \label{fica}  2^{-j'}  \geq j^2 2^{- j/ \eta } .  \EE
 We pick intervals of length $l = {j'}^2 2^{- \eta j'}$ and, for  the constant $A$ in Stute's lemma, we pick  $A = j$. 
Then (\ref{11bis})   applied with $N = [ 2^{\eta j'} ]$ yields that, with probability at  least $1- e^{-j^2} $,  the number of  intervals $\la'$ of scale $j'$  picked in such intervals is 
\[  2^{\eta j'} l  + O( j^2) =     O(j^2) . \]

We now assume that 
\BE \label{seca}   2^{-j'}  \leq j^2 2^{-j/\eta } . \EE
Then we pick
intervals of length $l = 3 \cdot 2^{-j} $, and $A = j+j' $. 
Then (\ref{11bis})   applied with $N = [ 2^{\eta j'} ]$ yields that, with probability at  least $1- e^{-(j+j')^2} $,  the number of  intervals $\la'$ of scale $j'$  picked in such  intervals is 
\BE \label{estimjp}  2^{\eta j'} l  + O( (j+j')^2 \sqrt{2^{\eta j'} l } ) \leq  2\cdot 2^{\eta j'} l. \EE

We are now ready to estimate the size of $l_\lambda$, assuming that all events described above happen 
(indeed,  we note that  the probabilities such that these events do not happen have a finite sum, so that, by the Borel-Cantelli lemma, they  a.s. all occur for $j$ large enough).

At scales $j'$ which satisfy (\ref{fica}), if at least one of the $\la ' \subset 3 \la$ does not vanish, then there are at most $j^2$ of them, and the corresponding contribution  to the sum in (\ref{eqllam}) lies between $| (c_{\la'})^p 2^{-(j-j')}| $ and $j^2 | (c_{\la'})^p 2^{-(j-j')}| $. 
At scales $j'$ which satisfy (\ref{seca}), the  contribution   of the wavelet coefficients of scale $j'$ to the sum  lies between $2^{\eta j' l} | (c_{\la'})^p 2^{-(j-j')}| $ and its double.  Since $c_{\la'} =  2^{-\al j'} $, the condition $ p < p_0$ implies that these quantities decay geometrically, so that the order of magnitude of the $p$-leader  is given by the first non-vanishing term in the sum.   Hence Proposition \ref{sizepleadlws}  holds.  \\

 \subsection{\boldmath$p$-exponents and lacunarity}

 We now derive the consequences of  Proposition \ref{sizepleadlws}  for the determination  of the $p$-exponents of  $X_{\al, \eta}$  at every point. 
 We first determine the range of $p$-exponents. First, note that all $p$-leaders have size at most $2^{- \al j}$, so that  the $p$-exponent is everywhere larger than $\al$.  In the opposite direction, as a  consequence of (\ref{estimjp}),  every interval $3 \la$ of scale $j$ includes at least one nonvanishing wavelet coefficient at scale $j/\eta + (\log j)^2$;
 therefore,  all $p$-leaders have size at least 
 \[ 2^{- \al \left( \frac{j}{\eta} + \log j \right)^2 - \frac{1}{p}\left( \frac{j}{\eta}  -j +( \log j)^2  \right)}. \]
  It follows that the $p$-exponents are  everywhere smaller than
\BE \label{Hmax}    H_{max} : = \frac{\al}{\eta} + \left(\frac{1}{\eta}  -1 \right) \frac{1}{p} . \EE
We have thus obtained that 
\[ \mbox{ a.s.}\qquad  \forall p < p_0, \; \forall x_0 \in \RR, \qquad  \al \leq h^p_{X_{\al, \eta}} (x_0) \leq H_{max}. \]

 For each $j$, let $E^j_\ome$ denote the subset of $[0,1] $ composed of intervals $3 \la$ ($\la \in \La_j$) inside which the first nonvanishing wavelet coefficient 
 is attained at a scale $l \leq [\ome j ]$, and let 
 \[ E_\ome = \limsup E^j_\ome. \]
 Proposition \ref{sizepleadlws} implies  that, if $x_0\notin  E_\ome $, then, for $j$ large enough, all wavelet leaders $l_{\la_j (x_0)}$ are bounded by 
   \[ j^2 2^{- \al  \frac{j}{\eta}  - \frac{1}{p}\left( \frac{j}{\eta}  -j \right)} , \] so that: 
   \BE \label{minopex}   \mbox{ if $x_0\notin  E_\ome $, then } \qquad h^p_{X_{\al, \eta}} (x_0) \geq {\al}{\ome} +  \frac{\ome -1}{p} . \EE
 On other hand, if  $x_0 \in  E_\ome $, then there exists an infinite number of $p$-leaders $l_{\la_j (x_0)}$ larger than 
 \[ 2^{- \al  \frac{j}{\eta}  - \frac{1}{p}\left( \frac{j}{\eta}  -j \right)} , \] so that: 
   \BE \label{majopex}  \mbox{ if $x_0\in  E_\ome $, then } \qquad h^p_{X_{\al, \eta}} (x_0) \leq {\al}{\ome} +  \frac{\ome -1}{p} . \EE
   It follows from (\ref{minopex}) and (\ref{majopex}) that   the sets of points where the $p$-exponent takes the value  
   \[ H = {\al}{\ome} +  \frac{\ome -1}{p} \] 
   are the sets 
   \[ H_\ome =  \bigcap_{ \omega' > \ome} E_{ \ome' }- \bigcup_{ \omega' < \ome} E_{ \ome' }. \] 
   We have thus obtained the following result.


 \BP   \label{propregloc} 
  Let $\al \in \RR$, $\eta \in (0,1)$ and 
Let $X_{\al, \eta}$ be a  lacunary wavelet series  of parameters $(\al, \eta)$. Let  $ p < p_0$;  the sets of points with a given $p$-exponent are the sets 
\[ F^p_{X_{\al, \eta}} (H) = H_\ome \quad \mbox{ for }  \quad \ome = \frac{H + 1/p}{\al + 1/p}; \] 
and additionally, if $x_0 \in H_\ome$, then 
\[ { \mathcal L}_{X_{\al, \eta}} (x_0)= \ome -1 . \] 
 \EP
 
 { \bf Remark:} 
 We actually do not need the wavelet used to be in the Schwartz class for Theorem \ref{theospec}  to be true. One can verify that, if the uniform regularity of the wavelet is larger than $\max ( |\al | , |H_{max}| )$, then all previous computations remain valid. \\

In order to determine the $p$-spectra and the lacunarity spectrum, one has to determine the Hausdorff dimensions of the sets $H_\ome$. We note that these sets 
 do not depend on $\al$ and on $p$, but only on the parameter $\ome$ and on the random drawing of the locations of the non-vanishing wavelet coefficients.  When $\al >0$, the dimensions of these sets  (expressed in a slightly different way) were determined in \cite{Jaf6}, where it is shown that   \[ \dim (H_\ome) = \eta \ome . \]
 The following result follows. 
 
 \BT \label{theospec}  Let $\al \in \RR$, $\eta \in (0,1)$ and 
let $X_{\al, \eta}$ be a  lacunary wavelet series  of parameters $(\al, \eta)$;  the  $p$-spectrum of $X_{\al, \eta}$ is supported by the interval $[\al, H_{max}]$ 
and, on this interval,
 \[ \mbox{ a.s. } \;\; \forall p < p_0, \;\; \forall H, \qquad d^p (H) = \eta \frac{H + 1/p}{\al + 1/p} .  \] 
Furthermore,  its lacunarity spectrum is given by 
 \[ \mbox{ a.s. } \;\; \forall  L\in[0, 1/\eta-1], \qquad d^\mathcal{L} (L) = \eta (L+1). \] 
 \ET
 
  { \bf Remark:} It is also shown  in \cite{Jaf6} that all the sets $H_\ome$ are everywhere dense, so that the quantity (\ref{oscil}) is equal everywhere  
  to $H_{max} -\al$. \\ 
 

For the sake of completeness, we now sketch how these dimensions can be computed. 
   We start by  estimating the size of $E_\ome$.  Note that  the number  of intervals $3 \la$ which comprise  $E^j_\ome$ is bounded by
   \[ [2^{\eta j} ] +  [2^{\eta (j+1)} ] + \cdots + [2^{\eta [\ome j] } ] \leq C 2^{\eta \ome j} . \]
   Using these intervals for $j \geq J$ as an $\ep$-covering, we obtain the following bound for the Hausdorff dimension of $E_\ome $
   \BE  \dim ( E_\ome ) \leq  \eta \ome . \EE
   We now consider the sets $ H_\ome $;  
   it follows from (\ref{minopex}) and (\ref{majopex}) that 
\[ H_\ome =  \bigcap_{ \omega' > \ome} E_{ \ome' }- \bigcup_{ \omega' < \ome} E_{ \ome' }. \] 
Since  $\forall \ome ' < \ome$, $ H_\ome \subset E_{\ome '} $, it follows from (\ref{majopex}) that  \[ \dim (H_\ome) \leq \eta \ome . \]

In order to get a lower bound on the Hausdorff dimension of $H_\ome $, we will need the following (slightly) modified notion of $\delta$-dimensional Hausdorff measure. 

\BD  \label{defmeshausbis} Let    $A\subset \RR$. 
For  $\varepsilon>0$  and  $\delta \in [0,1]$, let
\[ M^{\delta, \gamma}_{\varepsilon  } (A) =\inf_R \;  \left( \sum_{  i} | A_i |^\delta | \log (| A_i |)|^\gamma \right) ,\]
where $R$ denotes an  $\ep$-covering of    $A$, and where the infimum is taken on all 
 $\ep$-coverings.
The  
$(\delta, \gamma )$-dimensional Hausdorff measure of 
$A$ is 
\BE \label{defmesdel} M^{\delta, \gamma} (A) = \displaystyle\lim_{\ep\rightarrow 0}  M^{\delta, \gamma}_{\varepsilon, } (A)  . \EE
\ED

Since $E^j_\ome$  is composed of $\sim C 2^{\eta \ome j}$ randomly located intervals of length $3\cdot 2^{-j}$, 
 standard ubiquity arguments (such as in \cite{Jaf6, BerVel})   yield that 
 \[ M^{\eta \ome, 2}   (G_\ome )>0; \]  
(\ref{majopex}) implies that $\bigcup_{ \omega' < \ome} E_{ \ome' }$ (which  can be rewritten as a countable union)  has a vanishing $({\eta \ome, 2}) $-dimensional  Hausdorff measure. Thus 
\[ M^{\eta \ome, 2} \left( E_\ome - \bigcup_{ \omega' < \ome} E_{ \ome' }  \right) >0. \] 
Since this set is included in $H_\ome$, we obtain that 
 \[ \dim   (H_\ome )\geq\eta \ome .  \]

\begin{figure}[ht!!]
\centering
\includegraphics[width=\linewidth]{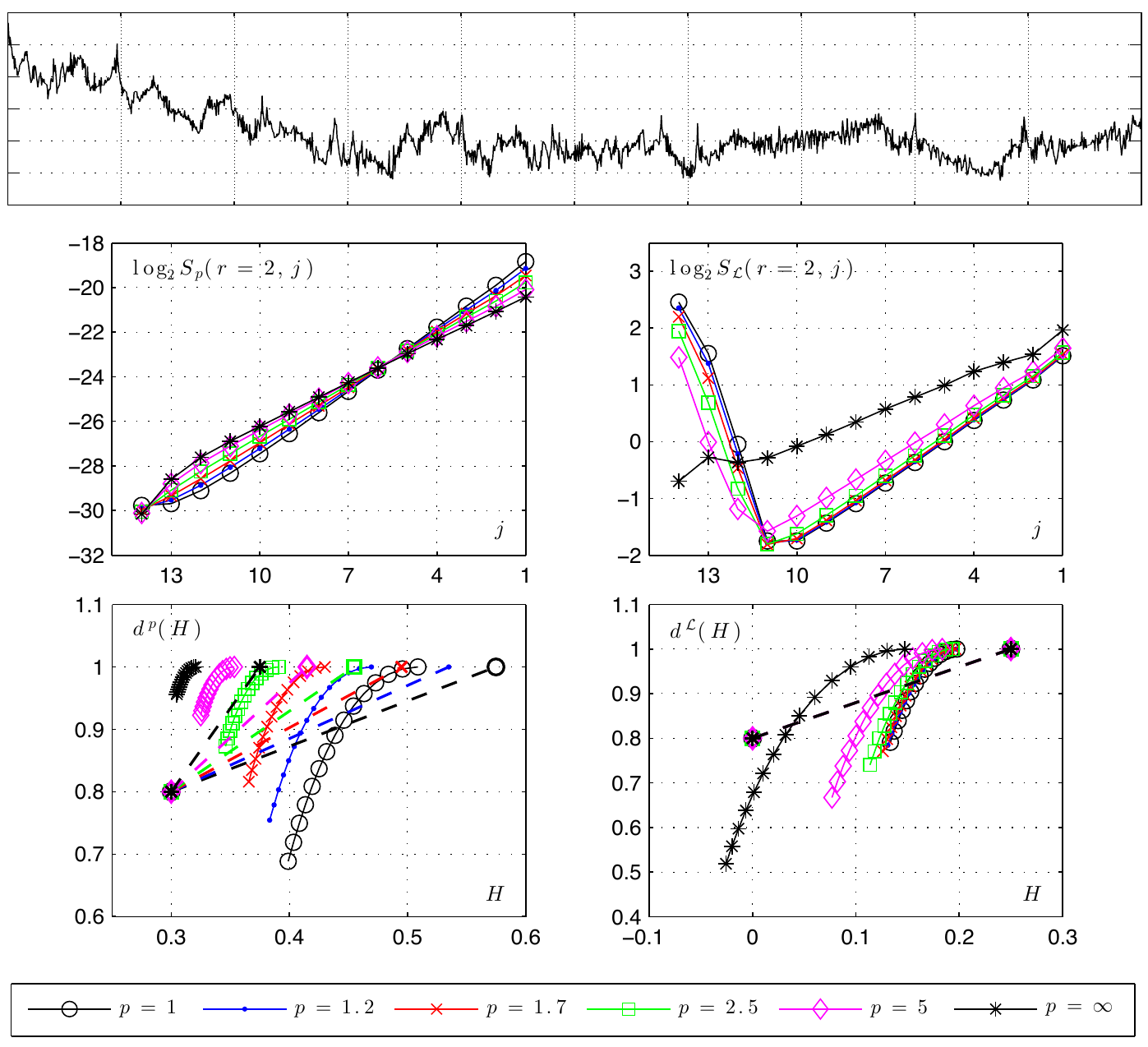}
\caption{\label{fig:LWS}{\bf Lacunary wavelet series}: A typical sample path of a lacunary wavelet series ($\alpha=0.3$, $\eta=0.8$, top row) and 
estimated structure functions (center row) and multifractal spectra (bottom row) for $p$-exponents (left column) and lacunarity exponents (right column) obtained with different values of $p$. The dashed lines indicate the theoretical multifractal spectra.}
\end{figure}

%

 
 It suffices now to rewrite these dimensions as a function of the $p$-exponent  to obtain  Theorem \ref{theospec}. \\

 Numerical examples for the estimation of $d^p (H)$ and $d^\mathcal{L} (H)$ of a lacunary wavelet series are given in Figure \ref{fig:LWS}. As predicted by theory, the numerical estimates of the $p$-exponent multifractal spectra are not invariant with $p$ but follow the evolution with $p$ of the theoretical spectra $d^p (H)$. The positions of the mode of the estimated spectra have a constant negative bias; yet, quantitatively, they  very well reproduce the shift of the mode of the theoretical spectra to smaller values of $H$ for increasing $p$, revealing the lacunary nature of the function. A refined analysis is possible with the estimated lacunarity exponent multifractal spectrum $d^\mathcal{L} (H)$, which has been computed here for several values of $p$ for illustration purposes. The mode of the spectrum is estimated at $H\approx 0.2$ (instead of the theoretical $H=0.25$). This clearly indicates the existence of positive lacunarity exponents. While the estimates for small values of $p$ fall short of revealing the full support of the theoretical multifractal spectrum, they still enable one to identify a relatively large interval of positive lacunarity exponent values. The best estimate of $d^\mathcal{L} (H)$ is obtained for the canonical value $p=p_0=+\infty$ ($q=q_0=0$) in this example and produces a satisfactory concave envelope of the theoretical multifractal spectrum that provides clear evidence for ensembles of lacunary singularities with a range of positive exponents.

\section{Concluding remarks}

The analysis that we developed is based on the assumption that   $p_0 (f)  >0$,  or that $\eta_f(p)>0$ for $p$ small enough, so that $p$-exponents can be defined, at least,  for $p \leq p_0$; we saw that this assumption allows us to deal with distributions of arbitrarily large order and, equivalently,  to model pointwise singularities with arbitrarily large  negative exponent. However, this does not imply that any tempered distribution satisfies these assumptions. 
Simple counterexamples  are supplied by the { \bf Gaussian fractional noises} $B_{\al}$  for $\al <0$ whose  sample paths can be seen as fractional derivatives of order $\frac{1}{2} - \al$ of the  sample paths of a Brownian motion on $\RR$ (Gaussian white noise corresponds to $\al =-1/2$, in which case it is a derivative, in the sense of distributions,  of Brownian motion).
In \cite{JAFFARD:2006:A}  the wavelet and leader scaling functions are derived, and it is proved that $\eta_{B_{\al}} = -\al p$, hence always is negative. 
However, the following result shows that, as soon as the wavelet expansion of the data has some sparsity, then this phenomenon no more occurs, and $p_0$ is always strictly positive (note that this situation is quite common in practice since sparse wavelet expansions are  often met in applications).  

\BD
A wavelet series $\sum_{j,k} c_{j,k} \psi_{j,k} $ is sparse if there exist $C >0$ and $\eta <1$ such that, on any interval $[l, l+1]$, 
\[ Card \{ k: \; c_{j,k} \neq 0  \} \leq C 2^{\eta j}. \]
\ED

Typical examples of sparse wavelet series are  supplied by lacunary wavelet series or by the measures which satisfy (\ref{hypsupp}). The following proposition implies  that multifractal analysis based on $p$-exponents is always possible for data with a sparse wavelet expansion. 

\BP \label{propspars} Let $f$ be a tempered distribution, which has a sparse wavelet expansion, then
$ \eta_f (p) >0 $  for $p $  small enough, so that $  p_0 (f) >0. $
\EP

 { \bf Proof:}  Since $f$ is a tempered distribution, it has a finite order, and thus it is a derivative of order $A$ of a continuous function. Therefore $f$ belongs to $C^{-A} (\RR)$, so that
  \[ |  c_{j,k} | \leq C 2^{Aj} . \]
   Using again compactly supported wavelets, the same argument as in the proof of Proposition \ref{propregmes}  yields that there are at most $C  2^{\eta j}$ nonvanishing wavelet coefficients at scale $j$; it follows that 
   \[ 2^{-j } \sum_{k} |  c_{j,k} |^p \leq C 2^{-j }  2^{\eta j }  2^{Apj }  \] 
 so that $\eta_f (p)  \geq 1-\eta -Ap$, and $\eta_f(p)  > 0$  for $p < (1-\eta )/A$.

\bibliographystyle{plain}
\bibliography{Lacun}

\end{document}